%Authors: M. Cepedello Boiso

%Title: Approximation of Lipschitz functions by $\Delta$-convex
%functions in Banach spaces

%Filename: boisoapproxlipsconv.tex
%TeX: AMSTeX
%Length: 43570 bytes
%Received Date: 2/13/97
%SubjectClass: 46B20
%Abstract: In this paper we give some results about the approximation 
%of a Lipschitz function on a Banach space by means of 
%$\Delta$-convex functions. In particular, we prove  that 
%the density of $\Delta$-convex functions in the set of 
%Lipschitz functions for the topology of uniform convergence
%on bounded sets characterizes the superreflexivity of the 
%Banach space. We also show that Lipschitz functions on 
%superreflexive Banach spaces are uniform limits on the 
%whole space of $\Delta$-convex functions.

%Citation: Preprint

%Special character check block
%32   space        33 ! exclam. pt.   34 " double quote  35 # sharp
%36 $ dollar       37 % percent       38 & ampersand     39 ' prime
%40 ( left paren.  41 ) rt. paren.    42 * asterisk      43 + plus
%44 , comma        45 - minus 46 . period        47 / divide
%58 : colon        59 ; semi-colon    60 < less than     61 = equal
%62 > greater than 63 ? question mark 64 @ at
%91 [ left bracket 92 \ backslash     93 ] right bracket 94 ^ caret
%95 _ underline    96 ` left single quote
%123 { left brace  124 | vertical bar 125 } right brace  126 ~ tilda

%

%----------------------------------------------------------------------------
% 
%----------------------------------------------------------------------------
%
\input amstex
\documentstyle{amsppt}
%
%
%
%
%
%
%\NoBlackBoxes
\NoRunningHeads
\TagsOnRight
\magnification=1200
%\loadbold
%%%%%%%%%%%%%%%%%%%%%%%%%%%%%%%%%%%%%%%%%%%%%%%%%%%%%%%%%%%%%%%%%%%%%%%%%%%%
%%%%%%%%%% TeX macros.
%
\def\roundhat{\!\raise0.333333333eX\hbox{$\smallfrown$}\!}
\def\subminusone{{\raise0.03333333ex\hbox{-}}\!1}
%
%%%%%%%%%%%%%%%%%%%%%%%%%%%%%%%%%%%%%%%%%%%%%%%%%%%%%%%%%%%%%%%%%%%%%%%%%%%%%
%%%%%%%%%% Title, abstract,...
\topmatter

\title 
Approximation of Lipschitz functions by $\Delta$-convex functions 
in Banach spaces
\\
\endtitle

\author 
Manuel Cepedello Boiso 
\endauthor

\address
Equipe d'Analyse, Universit\'e Pierre et Marie Curie--Paris 6, Paris.
\newline
\indent
Departamento de  An\'alisis Matem\'atico, Universidad de Sevilla, Sevilla.
\endaddress
\curraddr
Department of Mathematics, University of Missouri-Columbia, Columbia.
\endcurraddr

\email
%\newline 
%\indent 
manuel\@lebesgue.math.missouri.edu, cepedel\@ccr.jussieu.fr,  
boiso\@cica.es
\endemail

%  Math Subject Classifications 

\abstract 
In this paper we give some results about the approximation of 
a Lipschitz function on a Banach space
by means of $\Delta$-convex functions. In particular, we prove 
that the density
of $\Delta$-convex functions in the set of Lipschitz functions 
for the topology
of uniform convergence on bounded sets characterizes
the superreflexivity of the Banach space.
We also show that Lipschitz functions on superreflexive Banach 
spaces are uniform limits on the whole space of $\Delta$-convex 
functions.
\endabstract

\date
February 13, 1997
\enddate

\thanks
The author was supported by a FPU Grant of the Spanish {\it Ministerio de 
Educaci\'on y Ciencia}.   
\endthanks

\endtopmatter
%%%%%%%%%%%%%%%%%%%%%%%%%%%%%%%%%%%%%%%%%%%%%%%%%%%%%%%%%%%%%%%%%%%%%%%%%%%%
%%%%%%%% Text
\document

%%%%%%%%
\head 0.
Introduction and Notations
\endhead
%%%%%%%%

A function defined on a Banach space $X$ is called {\sl $\Delta$-convex} if 
it can be expressed as a difference of continuous convex functions or, 
equivalently,
if it belongs to the linear span of the continuous convex functions on $X$. 
The purpose of this paper is to give some necessary and sufficient 
conditions for the  approximation of Lipschitz functions by
$\Delta$-convex functions. 
The initial motivation for our work comes from two recent
articles of R. Deville, V. Fonf and P. H\'ajek (\cite {DFH$_1$} and 
\cite {DFH$_2$}). A consequence of their results is that, under
certain conditions on the Banach space $X$, any convex function on $X$
which is bounded on bounded sets can be approximated by smooth convex 
functions. It is therefore natural to consider 
the class of Banach spaces for which the $\Delta$-convex functions 
are dense in the class of Lipschitz functions
in order to extend this property of smooth approximations.
\par
Our main result is the
following characterization of superreflexivity.
\par
\proclaim{Theorem 0}
Let $X$ be a Banach space. Then $X$ is superreflexive if and only if 
every Lipschitz function on $X$ can be approximated uniformly on bounded
sets by differences of convex functions on $X$ which are bounded on bounded 
sets.
\endproclaim
A consequence of {\bf Theorem 0} is that the above mentioned approach 
does not provide any new result,
because it works only for superreflexive spaces (for which the 
smooth approximation property is known using the existence of 
partitions of unity, see Ch. VIII of \cite {DGZ}; 
for the analytic approximation case, see \cite {K}). 
\par%
For the superreflexive case, we give explicit formulas for the 
$\Delta$-convex 
approximation  of a Lipschitz function. These formulas, which are 
simpler than those from \cite {St}, also provide uniform convergence on 
the whole space $X$.
Specifically, the degree of convergence given by our formulas relies 
directly on the rotundity of the equivalent norm that it is used. 
Similar ideas in this
direction can be found in \cite {A} and \cite {PVZ}.
\par
We thank P. H\'ajek for bringing to our attention the
link between our work and the {\sl distortion theorem} (see \cite {OS} for
definitions and details). 
Our results provide simple formulas for deducing, on minimal superreflexive
Banach spaces (such as $\ell_p$, $1<p<\infty$), the existence of a convex
function which is not oscillation stable from the existence of a Lipschitz
function which is not oscillation stable.
\par
Let us fix some notation used in this paper. For a  real 
Banach space $X$, we denote an equivalent norm on $X$ by $\|\cdot\|$ 
and by $B_X$ its closed unit ball under this norm. 
By convex function we will always mean
continuous convex function. 
We will consider two fundamental topologies on the set of 
continuous functions defined on $X$: $\tau_\kappa$ (respectively $\tau_b$) 
is the topology of uniform  convergence on compact sets of $X$
(resp. uniform convergence on bounded sets of $X$).
\par
The {\sl modulus of convexity} of the norm $\|\cdot\|$ defined  by
$$
\delta_{\|\cdot\|}(\varepsilon)
= \inf
\bigg\{1-\bigg\|\frac{x+y}{2}\bigg\|:x,y\in B_X;\|x-y\|\geq\varepsilon\bigg\}
\quad 
\text {($0<\varepsilon<2$)}
$$
is called {\sl of power type $p$} ($p\geq 2$) if
$
\delta_{\|\cdot\|}(\varepsilon)\geq K\varepsilon^p$, for some $K>0$. 
\par
The concept of {\sl dyadic tree} will play an important role in the 
second part of this work.
Our trees are geometric trees contained in $X$ and defined as follows.
The symbol $\alpha$ denote a multi-index 
$\alpha=(\alpha_1\roundhat\alpha_2\roundhat\dots\roundhat\alpha_n)\in
{\{-1,1\}}^{<\Bbb N}$ and $|\alpha|:=n$.
For $n\in \Bbb N$, 
a dyadic $(n,\theta)$-tree $T$ in $X$ is a set of the form 
$
\big\{ x_{\alpha}\in X : \alpha \in {\{-1,1\}}^{<n} \big\}
$
satisfying the following two conditions:
\roster
\item  $x_{\alpha} = \frac {1}{2} x_{\alpha\roundhat 1}
                     + \frac {1}{2} x_{\alpha\roundhat\subminusone}$,
                     for all $|\alpha|<n$.
\item  $\|x_{\alpha} - x_{\alpha'}\| \geq \theta>0$, 
                       for $\alpha\neq\alpha'$.
\endroster
The point  $x_{\varnothing}$ will be called the {\sl root} of the tree $T$.
\par
%
%%%%%%%%
\head 1. 
The positive results
\endhead
%%%%%%%%

%
\proclaim{Theorem 1} 

Let $\roman{(}X,\| \cdot \| \roman{)}$ be a Banach space. The norm 
$\| \cdot \|$ is locally uniformly convex (respectively uniformly convex)
if and only if the following property holds:
for every Lipschitz function $f$ on $X$, the sequence of functions 
${(f_n)}_{n\in\Bbb N}$ defined by the formula
$$
f_n(x) := \inf_{y\in X} 
\Big\{ 
f(y) + n \big( {2\| x \|}^2 + 2 {\| y \|}^2 - {\| x+y \|}^2 \big) 
\Big\}  
\quad ( n\in \Bbb N,\  x\in X )
$$
is $\tau_\kappa$-converging (resp. $\tau_b$-converging) to $f$.
\endproclaim
\remark{Remark}
For any function $f$ on $X$ and $n\in \Bbb N$,  the function $f_n$ defined 
as above is a $\Delta$-convex function. This follows immediately from the 
decomposition  
$f_n = c_n - d_n$, with  
$
c_n(x) := 2n {\| x \|}^2 
$
%\quad \text {and} \quad
and
$$
d_n(x) := \sup_{y\in X} 
\big\{ n{\| x+y \|}^2 - 2n{\| y \|}^2 - f(y) \big\}
\quad (x \in X).
$$
The functions $c_n$ and $d_n$  are clearly convex.
\endremark
\demo{Proof of the Theorem 1}

Let us see first that the property is necessary. So, let $f$ be a Lipschitz 
function on $X$. We have to show that the previously defined sequence 
$(f_n)_{n\in \Bbb N}$ ($\tau_K$ or $\tau_b$)-converges to $f$ 
if the norm $\| \cdot \|$ satisfies the corresponding rotundity condition. 
\par 
We begin with the following general result:
\remark{Fact}
For any point $x\in X$, $(f_n(x))_{n\in \Bbb N}$ is an increasing sequence
bounded above by $f(x)$.
\endremark
This {\sl Fact} follows immediately from taking $y=x$ in the $f_n$'s
infimum formula and from the inequality 
$$
2 {\| x \|}^2 + 2 {\| y \|}^2 - {\| x+y \|}^2
\geq
{2\| x \|}^2 + 2 {\| y \|}^2 - {\big (\| x\| + \| y \|\big)}^2 
=
{\big( \| x \| - \| y \| \big)}^2
\geq
0.
\tag{1}
$$
\par
Since for $K \in \Bbb N$ we have that 
$ 
f_{(K\cdot n)}
=
K \big( \frac {f}{K} \big)_n
$,
the previous {\sl Fact} allows us 
to suppose without loss of generality that the Lipschitz constant of
$f$ is less than $1$
({\it i.e.} $|f(x)-f(y)| \leq \| x - y \|$ for all $x,y \in X$).
\par
We need to  study of the infimum formula defining $f_n$ 
at a point $x\in X$.
Thus, consider any point $y$ for which we have 
$$
f(y) + n \big( {2\| x \|}^2 + 2 {\| y \|}^2 - {\| x+y \|}^2 \big) 
\leq 
f(x).
\tag {2}
$$ 
As $f$ is $1$-Lipschitz, we deduce from {\bf (1)} and {\bf (2)} that
$$
n {\big( \| x \| - \| y \| \big)}^2
\leq
n \big( 2{\| x \|}^2 + 2{\| y \|}^2 - {\| x+y \|}^2 \big)
\leq 
 f(x)-f(y) 
\leq
\| x - y \|.
\tag{3}
$$ 
\par
This last condition  gives a relation between the norms of $x$ and 
$y$. Actually, $\| y\|$ can be controlled by $\|x\|$ in the following way.
Suppose that $\| y\|  \geq 1 + \|x\|$, then by {\bf (3)}
we get that
$$
1
\leq
\big| \|x\| -\|y\| \big| 
\leq 
{\big (\| x\| - \| y \|\big)}^2
\leq 
\frac 1{n}\| x - y \|
\leq
\frac 1{n} \|x\| + \frac 1{n} \|y\|.
\tag{4}
$$
And we conclude for $n\geq 3$ that 
$
\| y\| \leq  \frac {n+1}{n-1}\|x\| \leq 2 \|x\|.
$
Therefore, we see that if $y$ satisfies {\bf(2)} and $n\geq 3$ then 
$$
\|y\| \leq 2(1+\|x\|).
\tag{5}
$$
Then, for $n\geq 3$ we deduce that 
$$
f_n(x) = \inf_{ \|y\|\leq 2(1+\|x\|) } 
\big\{ 
f(y) + {2n\| x \|}^2 + 2n {\| y \|}^2 - n{\| x+y \|}^2 
\big\}.
\tag{6}  
$$
In particular, the boundedness on bounded sets of $f_n$ 
(and, consequently, of $d_n=c_n - f_n$) follows immediately. 
\par
Moreover, the upper bound on $\|y\|$ given by {\bf (5)} 
together with the condition {\bf (3)} 
gives
$$
0
\leq
2 {\| x \|}^2 + 2 {\| y \|}^2 - {\| x+y \|}^2 
\leq
\frac {1}{n}\| x - y \|
\leq
\frac {1}{n} \|x\| + \frac {1}{n} \|y\|
\leq
\frac {3}{n} (1 + \|x\|).
\tag{7}
$$
\par
The condition {\bf (7)} gives the crucial step of the proof.
In fact, if the norm $\|\cdot\|$ verifies one of the rotundity properties 
stated in theorem
then by {\bf (7)} $n$ can be chosen big enough  to 
necessarily enforce $y$ to be close to $x$. Therefore, $f_n(x)$
must be close to $f(x)$. Let us justify this assertion.
\par
Suppose that the sequence of functions $f_n$ is not compactly converging
to $f$.
Since the sequence of $\Delta$-convex functions ${(f_n)}_n$ is increasing
and $f$ is continuous, Dini's theorem tells us that 
$\tau_\kappa$-convergence of $f_n$ to $f$ is equivalent to the pointwise
convergence. 
Then, there exists a  point $x_0\in X$ such that ${(f_n(x_0))}_n$ does 
not converge to $f(x_0)$. 
\par
As ${(f_n(x_0))}_n$ is increasing, 
there exists some $\varepsilon_0>0$ such that for any $n\in \Bbb N$
we have $f_n(x_0) + \varepsilon_0 < f(x_0)$. By definition of $f_n$,
we can find a sequence of ${(y_n)}_n$ in $X$ so that
$$
f_n(x_0) + \varepsilon_0
\leq 
f(y_n) + {2n\| x_0 \|}^2 + 2n {\| y_n \|}^2 - n{\| x_0 + y_n \|}^2
+\varepsilon_0
\leq 
f(x_0).
\tag{8}
$$
Since $f$ is $1$-Lipschitz, we get from {\bf (8)} that 
$$
\|x_0 - y_n \| 
\geq 
f(x_0) - f(y_n) 
\geq 
\varepsilon_0 
>0.
\tag{9}
$$
Since $y_n$ verifies the condition {\bf (8)}, it follows from {\bf (7)}
 that we have
$$
0
\leq
2 {\| x_0 \|}^2 + 2 {\| y_n \|}^2 - {\| x_0+y_n \|}^2 
\leq
\frac {3}{n} (1 + \|x_0\|)
@>> {n\to \infty}>   0.
\tag{10}
$$
But {\bf (9)} and {\bf (10)} show that the norm $\|\cdot\|$ can not be 
locally uniformly
convex at $x_0$ (cf .Ch. II. Prop. 1.2 of \cite {DGZ}).
That proves the compact convergence of $f_n$ to $f$ if the norm 
$\| \cdot \|$ is locally uniformly convex.
\par
A simple proof of the uniform convergence of the sequence $f_n$ 
in the uniformly convex case follows the same lines.
We will give later a direct quantitative approach (see the proof 
of {\bf Theorem 3}).
If the sequence $f_n$ does not converge uniformly on a bounded set of
$X$, there is an $\varepsilon_0>0$ and a bounded sequence ${(x_n)}_n$
so that $f_n(x_n)+\varepsilon_0 < f(x_n)$, $n\in \Bbb N$. Then,
for each $x_n$ ($n\in \Bbb N)$ we can choose $y_n$ verifying
$$ 
f(y_n) + {2n\| x_n \|}^2 + 2n {\| y_n \|}^2 - n{\| x_n + y_n \|}^2
+\varepsilon_0
\leq 
f(x_n).
\tag{11}
$$
\par
Similar reasonings as before imply from {\bf (11)} that the next 
two statements are fulfiled:
$$
\gather
\|x_n - y_n \| 
\geq 
f(x_n) - f(y_n) 
\geq 
\varepsilon_0 
>0
\tag{12}
\\ 
0
\leq
2 {\| x_n \|}^2 + 2 {\| y_n \|}^2 - {\| x_n+y_n \|}^2 
\leq
\frac {3}{n} (1 + \|x_n\|)
@>> {n\to \infty}>   0
\tag{13}
\endgather
$$
\par
We have 
$
\lim_{n\to\infty} \frac {3}{n} (1 + \|x_n\|)=0$  since the sequence 
$(x_n)_n$ is bounded.
Moreover, {\bf (5)} implies that the sequence $(y_n)$ is also bounded.
Therefore, {\bf (12)} and {\bf (13)} show that the
norm $\|\cdot\|$ cannot be uniformly convex 
(cf. Ch. IV Lemma 1.5 \cite {DGZ}). 
The necessity of the property is proved. 
\par
Conversely, the property stated in the theorem is sufficient. We first 
show that the property of compact convergence implies 
the following claim.
\proclaim{Claim 1.1}
For any sequence ${(x_n)}_{n\in \Bbb N} \subset X$ and any point $x_0 \in X$,
the condition 
$
\lim_{n\to\infty}
\big( 2 {\| x_0 \|}^2 + 2 {\| x_n \|}^2 - {\| x_0 + x_n \|}^2 \big)
=
0
$
implies $\text{d}\big(x_0, (x_n)_n\big)=0$.
\endproclaim
Since the {\bf Claim 1.1} is also valid for any subsequence of the given
sequence $(x_n)_n$, we can strengthen the conclusion of the claim to 
$\lim_{n\to\infty} \|x_n-x_0\|=0$ and
the locally convexity of the norm $\| \cdot \|$ is verified
(cf .Ch. II. Prop. 1.2 of \cite {DGZ}).
\demo{Proof of the Claim 1.1}
Consider the Lipschitz function $g(x)=\text{d}\big(x,(x_n)_n\big)$, $x\in X$.
For $x_0\in X$ and $\varepsilon>0$, the pointwise convergence of the
sequence $(g_n)_n$ to $g$ at $x_0$ tells us
that there exists $N_\varepsilon\in \Bbb N$ so that 
$$
\text{d} \big( x_0,(x_n)_n \big)
=
g(x_0)
\leq \inf_{y\in X} 
\Big\{ 
g(y) + N_\varepsilon 
\big( 2 {\| x_0 \|}^2 + 2 {\| y \|}^2 - {\| x_0 + y \|}^2 \big)
\Big\}
+ \varepsilon.
\tag{14}
$$
\par
Using $g(x_n)=0$, we can evaluate 
the inequality {\bf (14)} at $y=x_n$ ($n\geq N_\varepsilon$) and deduce that
$$
\text{d} \big( x_0,(x_n)_n \big)
\leq 
N_\varepsilon 
\big( 
\lim_{n\to\infty} 2{\| x_0 \|}^2 + 2{\| x_n \|}^2 - {\| x_0 + x_n \|}^2
\big)
+ \varepsilon
=
\varepsilon.
$$
\par
Therefore, $\text{d}\big(x_0, (x_n)_n\big)=0$ 
and the claim is proved.
\enddemo
\par
If the property of uniform convergence on bounded sets of $X$ holds,  
the next claim, analogous to the previous one, also does. 
\proclaim{Claim 1.2} 
Let $(x_n)_{n\in \Bbb N}$ and $(y_n)_{n\in \Bbb N}$ be two bounded sequences
in $X$ such that 
$
\lim_{n\to\infty} 
\big( 2 {\| x_n \|}^2 + 2 {\| y_n \|}^2 - {\| x_n + y_n \|}^2 \big) = 0
$. Then $\lim_{m\to\infty} \text{d} \big( y_m, (x_n)_n \big) = 0$.
\endproclaim
\demo{Proof of the Claim 1.2}
As before, consider the Lipschitz function 
$\text{d} \big( \cdot, (x_n)_n \big)$ and $\varepsilon>0$.
This time, the property of uniform convergence gives a positive
integer $N_\varepsilon$ so that for all $m\in \Bbb N$ the next inequality
is satisfied.
$$
\text{d}\big( y_m, (x_n)_n \big)
 \leq 
\inf_{z\in X}
\Big\{ 
\text{d} \big( z,(x_n)_n \big) + N_\varepsilon 
\big( 2 {\| z \|}^2 + 2 {\| y_m \|}^2 - {\| z + y_m \|}^2 \big)
\Big\}
+ \varepsilon.
\tag{15}
$$
Taking $z=x_m$ in the infimum of {\bf (15)} we obtain
$$
\text{d}\big( y_m, (x_n)_n \big)
\leq
N_\varepsilon
\big( 2 {\| x_m \|}^2 + 2 {\| y_m \|}^2 - {\| x_m + y_m \|}^2 \big)
+ \varepsilon
\leq
2 \varepsilon,
$$
for $m$ large enough, and the claim is proved.
\enddemo
In order to finish with the proof, we shall show that the validity of
{\bf Claim 1.2} implies that the norm $\| \cdot \|$ is uniformly convex.
If not, (cf. Ch. IV Lemma 1.5 \cite {DGZ}) there exists two bounded 
sequences $(x_n)_{n\in \Bbb N}$ and $(y_n)_{n\in \Bbb N}$ in $X$ satisfying
that 
$$
\lim_{n\to\infty}
\big( 2 {\| x_n \|}^2 + 2 {\| y_n \|}^2 - {\| x_n + y_n \|}^2 \big)
= 0
\text{ and } 
\|x_n - y_n\| \geq 1 \text { (for all $n\in \Bbb N$).}
\tag{16}
$$
First, the norm $\|\cdot\|$ is locally uniformly convex (since 
{\bf Claim 1.1} clearly holds). It follows that the sequence $(x_n)_n$ 
has no norm cluster point. Indeed, if for some $x_0\in X$
there exists $(x_{n_k})_k @>>k> x_0$, then 
$$
\lim_{k\to\infty} 2\|x_0\|^2 + 2\|y_{n_k}\|^2 - \|x_0+y_{n_k}\|^2
=
\lim_{k\to\infty} 2\|x_0\|^2 + 2\|x_{n_k}\|^2 - \|x_0+x_{n_k}\|^2
= 0
$$
and therefore $(y_{n_k})_k @>>k> x_0$. A contradiction with the
fact $\|x_n - y_n\| \geq 1$, for all $n\in \Bbb N$, of {\bf (16)}.
\par
Hence, passing to a subsequence, we can suppose that for some $1>\alpha>0$
we have that $\|x_n - x_m\| \geq \alpha$, for all $n\not= m$.
Now, we need the following technical lemma, whose proof will be given later.
\proclaim{Lemma 1.3}
Let be $(\lambda_n)_{n\in \Bbb N}\subset [0,1]$ and 
$(x_n)_{n\in \Bbb N}$, $(y_n)_{n\in \Bbb N}$ two bounded sequences
in $X$ so that
$
2\|x_n\|^2 + 2\|y_n\|^2 - \|x_n + y_n\|^2 @>>n\to\infty> 0
$. 
Then the sequence of convex combinations
$
z_n=\lambda_n x_n + (1-\lambda_n)y_n
$ 
($n\in \Bbb N$) satisfies that
$
2\|x_n\|^2 + 2\|z_n\|^2 - \|x_n + z_n\|^2 @>>n\to\infty> 0
$.
\endproclaim
For any $n\in \Bbb N$, take $0 \leq \lambda_n \leq 1$ such that for
$z_n:=\lambda_n x_n + (1-\lambda_n) y_n$ we have 
$\|x_n - z_n\|=\frac {\alpha}{2}$. Then, the {\bf Lemma 1.3} and the
{\bf Claim 1.2} used jointly imply that 
$\lim_{m\to\infty} \text{d} (z_m,(x_n))=0$. But, 
$\|x_n - z_n\|=\frac {\alpha}{2}>0$ and also for $n\not= m$ 
$\|x_n - z_m\|\geq \|x_n - x_m\| - \|x_m - z_m\| \geq \alpha - 
\frac {\alpha}{2} = \frac {\alpha}{2}>0$, a contradiction.
\qed
\enddemo
\demo{Proof of the Lemma 1.3}
By the inequality
$$
0 
\leq 
\big(\|x_n\|-\|y_n\|\big)^2 
\leq  
2\|x_n\|^2 + 2\|y_n\|^2 - \|x_n + y_n\|^2
@>>n\to\infty> 0,
\tag{17}
$$
we have that
$$
\lim_{n\to\infty} \big( \|x_n\| - \|y_n\| \big)=0.
$$
Since the sequences $(x_n)_n$ and $(y_n)_n$ are bounded we also
have that 
$$
\lim_{n\to\infty} \big( \|x_n\|^2 - \|y_n\|^2\big) = 0.
\tag{18}
$$
But then we deduce from {\bf (17)} and {\bf (18)} that 
$$
\lim_{n\to\infty} 
\bigg( \|x_n\|^2 - {\Big\| \frac {x_n+y_n}{2} \Big\|}^2 \bigg)
=
\lim_{n\to\infty} 
\bigg( \|y_n\|^2 - {\Big\| \frac {x_n+y_n}{2} \Big\|}^2 \bigg)
= 0.
\tag{19}
$$
\par
On the other hand, using the convexity of the function ${\|\cdot\|}^2$
we get the following general lower estimate for $0\leq\lambda\leq 1$ 
$$
\| \lambda x + (1-\lambda) y \|^2
\geq
\bigg\| \frac {x+y}{2} \bigg\|^2\!\!-\! \big| 1 - 2\lambda \big| 
\max
\Bigg\{ 
\Bigg| \bigg\| \frac {x+y}{2} \bigg\|^2 - \|x\|^2 \Bigg| ,
\Bigg| \bigg\| \frac {x+y}{2} \bigg\|^2 - \|y\|^2 \Bigg|
\Bigg\}.
$$
\par
Putting {\bf (17)}, {\bf (18)}, {\bf (19)} and the last inequality together, 
we obtain that
$$
\split
0
&\leq 
2\|x_n\|^2 + 2\|z_n\|^2 - \|x_n+z_n\|^2
\\
&=
2\|x_n\|^2 + 2\|\lambda_n x_n + (1-\lambda_n)y_n\|^2 - 
\| (1+\lambda_n)x_n + (1-\lambda_n)y_n \|^2
\\
&\leq
2\|x_n\|^2 + 2\lambda_n\|x_n\|^2 + 2(1-\lambda_n)\|y_n\|^2 -
4 \bigg\| \frac {1+\lambda_n}{2}x_n + \frac {1-\lambda_n}{2}y_n \bigg\|^2
\\
&\leq
2\|x_n\|^2 + 2\|y_n\|^2 - \|x_n+y_n\|^2 +
\lambda_n \Big( \|x_n\|^2 -\|y_n\|^2 \Big) 
\\
&\quad + 4 \lambda_n
\max
\Bigg\{ 
\Bigg| \bigg\| \frac {x_n+y_n}{2} \bigg\|^2 - \|x_n\|^2 \Bigg| ,
\Bigg| \bigg\| \frac {x_n+y_n}{2} \bigg\|^2 - \|y_n\|^2 \Bigg|
\Bigg\}  
@>>n\to\infty> 0.
\endsplit
$$
The lemma is proved and this concludes the proof of {\bf Theorem 1}.
\enddemo
\remark{Remark}
The infimum formula used for the definition of $(f_n)$ 
in the {\bf Theorem 1} is closely related to the well-known 
{\sl inf-convolution formula} of $f$ by $n{\| \cdot \|}^2$:
$$ 
\big( f\phantom{[} \square\phantom{[} (n{\|\cdot\|}^2) \big) (x):=
\inf_{y\in X} \big\{ f(y) + n{\|x-y\|}^2 \big\}
\quad
\text { ($n\in \Bbb N$, $x\in X$)}.
\tag{20}
$$
In fact, these two infimum formulas are identical if the norm $\|\cdot\|$
is a Hilbertian norm,
because of the parallelogram identity. 
However, for a non-Hilbertian norm $\|\cdot\|$ the functions given by the 
inf-convolution formula
can not be expressed in general as $\Delta$-convex functions.
\endremark
\par
As a corollary of the above remark, we have that the formula of 
{\bf Theorem 1} converges uniformly on $X$ for case of a Hilbertian norm
$\|\cdot\|$ (since it is well-known that the 
inf-convolution formula of {\bf (20)} 
converges uniformly on the whole space $X$, see \cite {LL}).
The question that
naturally arises is whether this remains true or not for a general
uniformly convex norm. The answer is given in the following proposition.
\par
\proclaim{Proposition 2}
Let $(X,\|\cdot\|)$ be a Banach space. If 
for every Lipschitz function $f$ on $X$ the sequence of functions 
$$
f_n(x) := \inf_{y\in X} 
\Big\{ 
f(y) + n \big( {2\| x \|}^2 + 2 {\| y \|}^2 - {\| x+y \|}^2 \big) 
\Big\}  
$$
converges to $f$ uniformly on $X$, then the modulus of convexity of the
norm $\|\cdot\|$ is of power type 2.
\endproclaim
\demo{Proof of the Proposition 2}
This uniform convergence property has the following consequence analogous
to {\bf Claim 1.1} and {\bf Claim 1.2} and whose proof is identical to
that of {\bf Claim 1.2}.
\proclaim{Claim 2.1}
If $(x_n)_{n\in \Bbb N}$ and $(y_n)_{n\in \Bbb N}$ are two sequences
(not necessarily bounded) in $X$ satisfying that
$
\lim_{n\to\infty} 
\big( 2 {\| x_n \|}^2 + 2 {\| y_n \|}^2 - {\| x_n + y_n \|}^2 \big) = 0
$, then they also verify that the 
$\lim_{m\to\infty} \text{d} \big( y_m, (x_n)_n \big) = 0$.
\endproclaim
By {\bf Claim 1.2}, any norm which satisfies the conclusion of 
{\bf Claim 2.1} is uniformly convex. In fact, {\bf Claim 2.1} 
insures that the modulus of
convexity of the norm $\|\cdot\|$ is of power type 2. 
If not (see \cite {H}), there exists two sequences 
$(x_n)_{n\in\Bbb N}$ and $(y_n)_{n\in\Bbb N}$ such that
$$
2\|x_n\|^2 + 2\|y_n\|^2  < \|x_n+y_n\|^2 + \frac {1}{n} \|x_n-y_n\|^2
\quad 
\text {$x_n \not= y_n$, $n\in \Bbb N$}.
\tag {21}
$$
If we take $u_n:= \frac {x_n}{\|x_n-y_n\|}$ and 
$v_n:= \frac {y_n}{\|x_n-y_n\|}$ in {\bf (21)} we obtain that 
$$
\gather
0\leq 2\|u_n\|^2 + 2\|v_n\|^2 - \|u_n+v_n\|^2< \frac {1}{n} @>>n\to\infty>0,
\tag{22}
\\ 
\|u_n - v_n\| = 1.
\tag {23}
\endgather
$$
\par
Since the norm $\|\cdot\|$ is uniformly convex and conditions 
{\bf (22)} and {\bf (23)} holds, the sequences $(u_n)_n$ and
$(v_n)_n$ can not be bounded. Notice that we have also from {\bf (22)} that 
$\lim_{n\to\infty} \big( \|u_n\|-\|v_n\|\big)=0$. 
Thus, passing to a subsequence we can suppose that 
$$
\gather
\sup\{\|u_n\|,\|v_n\|\} +1 <\min\{\|u_{n+1}\|,\|v_{n+1}\|\}
\text { for all $n\in\Bbb N$ }
\\
\Longrightarrow
\|u_n - v_m\| \geq \big| \|u_n||-\|v_m\| \big| \geq 1
\quad
\text{($n\not=m$).}
\tag {24}
\endgather
$$
Assembling {\bf (23)} and {\bf (24)} we deduce that
$
\text{d}\big( v_m,(u_n)_n \big)=1 
$
($m\in\Bbb N$). 
But {\bf Claim 2.1} and {\bf (22)} imply that
$\lim_m\text{d}\big( v_m,(u_n)_n \big)=0$ which is a contradiction.
\qed 
\enddemo
\par
{\bf Theorem 3} below provides a converse to {\bf Proposition 2} and
an explicit formula for uniform approximation of Lipschitz functions by
$\Delta$-convex functions on superreflexive Banach spaces. 
The validity of this
formula depends upon the existence in the superreflexive space of an
equivalent norm which is enough rotund.
\proclaim{Theorem 3}
Let $(X,\|\cdot\|)$ be a Banach space whose norm $\| \cdot \|$ has its
modulus of convexity of power type $p$ ($p\geq 2$). Then for every Lipschitz
function $f$ on $X$ the following sequence of $\Delta$-convex functions
$$
f_n^p (x)
=
\inf_{y\in X}
\Big\{ 
f(y) + n \big( {2^{p-1}\| x \|}^p + 2^{p-1}{\| y \|}^p - {\| x+y \|}^p \big) 
\Big\}  
\quad ( n\in \Bbb N,\  x\in X )
$$
converges to $f$ uniformly on $X$.
\endproclaim
\demo{Proof of the Theorem 3}
We essentially need the following lemma to prove the result.
\proclaim{Lemma 3.1}
If the  modulus of convexity of the norm $\|\cdot\|$ is of power type $p$
($p\geq 2$) then there exists a positive constant $C_{\|\cdot\|}\leq 1$ 
such that for every pair $x,y \in X$ the following inequality holds:
$$
C_{\|\cdot\|} {\| x - y \|}^p
\leq
 {2^{p-1}\| x \|}^p + 2^{p-1}{\| y \|}^p - {\| x+y \|}^p .
$$
\endproclaim
\demo{Proof of the Lemma 3.1}
Under the assumption of the lemma (see \cite {H}), 
there exists a positive
constant $C'_{\|\cdot\|}\leq 2$ such that for every pair
$u,v\in X$ we have that
$$
\|u+v\|^p + \|u-v\|^p \geq 2\|u\|^p +C'_{\|\cdot\|}\|v\|^p.
\tag{25}
$$ 
The lemma follows from the change of variables
$x=\frac {u+v}{2}$ and $y=\frac {u-v}{2}$ in {\bf (25)}.
\qed
\enddemo
Given a Lipschitz function $f$ on $X$, we use 
the previous {\bf Lemma 3.1} and the fact that
$f_n^p(x)\leq f(x)$ (for all $x\in X$ and $n\in\Bbb N$) to obtain the
following chain of inequalities:
$$
\big( f\phantom{[} \square\phantom{[}(nC_{\|\cdot\|} {\|\cdot\|}^p) \big)
\leq
f_n^p
\leq
f
\quad
(n\in \Bbb {N}).
\tag{26}
$$
Since $p\geq 2$, the previous inf-convolution formula in  {\bf (26)}
converges to $f$ uniformly on $X$ (see \cite {LL}). 
Therefore, the same is true for $(f_n^p)_n$.
For instance, if $f$ is a $1$-Lipschitz function we get the following 
estimate
$$
\|f-f_n^p\|_{\infty}
\leq\
\|f-
\big(f\phantom{[}\square\phantom{[}(nC_{\|\cdot\|}{\|\cdot\|}^p)\big)
\|_{\infty}
\leq
\bigg ( \frac {1} {nC_{\|\cdot\|}}\bigg)^{\frac {1}{p-1}}.
\qed
$$
\enddemo
\par
%
%%%%%%%%%%%%%%%%%%%%%%%%%%%%%%%%%%%%%%%%%%%%%%%%%%%%%%%%%%%%%
%
The first consequence that stems from the {\bf Theorem 1} is another
proof of a result due to G.A. Edgar (\cite {E}).
\par
\proclaim{Corollary 4}
Let $X$ Banach space with an equivalent locally uniformly convex norm. 
Then the $\sigma$-fields of Borel sets for the norm and weak 
topologies are the same. 
\endproclaim
\demo{Proof of the Corollary 4}
The first part of the {\bf Theorem 1} affirms that  the existence
of a locally uniformly convex  norm on $X$ implies the following:
$$
\forall f:X \to \Bbb R \text {\ Lipschitz\ } \exists {(c_n)}_n, {(d_n)}_n
\subset \text{Conv}(X) \text {\ such that\ }
f=\tau_\kappa \text{-} \lim_{n\to \infty} (c_n - d_n).
$$
Let us check that this property implies the equivalence of the two Borel
families, $\text{Bor}(X, \| \cdot \|)$ and
$\text{Bor}(X, w )$. Obviously,
$
\text{Bor}(X, w)
\subseteq
\text{Bor}(X, \| \cdot \| )
$. 
To see the other inclusion, take $F$ a $\|\cdot\|$-closed set of $X$.
Then consider the Lipschitz function $f(\cdot)=\text{dist}( \cdot, F)$.
Note that the previous property implies that $f$ is the pointwise limit
of a sequence of $w$-Borel functions 
(because every continuous convex function
is $w$-lower-semicontinuous). Therefore, $f$ is $w$-Borel and so
$F$ is $w$-Borel. 
\qed
\enddemo
\par
Now, we proceed to state two applications of 
{\bf Theorem 1} and {\bf Theorem 3} to the study of 
superreflexive Banach spaces.
Actually, we will show in the next section that both of them 
characterize superreflexivity.
\proclaim{Corollary 5}
Let $X$ be a superreflexive Banach space. Denote by 
${\roman{Conv}} \roman{(} X \roman{)}$ the set of continuous convex 
functions on $X$, ${\roman{Conv}}_b \roman{(} X \roman{)}$ the subset of
${ \roman{Conv} \roman{(} X \roman{)} } $ consisting of the functions 
which are bounded on bounded sets and ${\Cal{UC}}_b(X)$ the class of functions
on $X$ which are uniformly continuous on bounded sets of $X$. Then 
$$
{\overline {\roman{span}}}^{\tau_b} 
\big\{ {\roman{Conv}}_b(X) \big\}
=
{\Cal{UC}}_b(X).
$$
\endproclaim  
\demo{Proof of the Corollary 5}
From the Prop.~1.6 of \cite {Ph} follows immediately that
$
{\text{Conv}}_b(X)
=
{\text{Conv}}(X) \cap {\Cal{UC}}_b(X) 
\subseteq 
{\Cal{UC}}_b(X)
$.
Let us show the $\tau_b$-density of the $\Delta$-convex functions in the
set ${\Cal{UC}}_b(X)$.
\par
Notice that the set of Lipschitz functions on $X$ is dense in 
${\Cal{UC}}_b(X)$ under the topology $\tau_b$. This fact can be proved using
again the inf-convolution formula.
More precisely,
given $f \in {\Cal{UC}}_b(X)$ and  bounded on $X$,
$ 
\big( f\phantom{[} \square\phantom{]} (n{\|\cdot\|}) \big) (x) 
$ 
is a sequence of Lipschitz functions $\tau_b$-converging  
to $f$.
Since the bounded functions of ${\Cal{UC}}_b(X)$ are clearly $\tau_b$-dense
in ${\Cal{UC}}_b(X)$,
the density of the Lipschitz functions 
is therefore deduced. 
\par
Thus,
the corollary holds if we show that the convex functions
$\{ c_n, d_n  \}_{n\in \Bbb N}$ 
obtained in the proof of {\bf Theorem 1}
are in ${\text{Conv}}_b(X)$. Clearly, 
$c_n(\cdot) = 2n {\| \cdot \|}^2 \in {\text{Conv}}_b(X)$
and, as we remarked during the proof of this theorem, $d_n=c_n - f_n$ is also
bounded on bounded sets of $X$.
\qed
\enddemo
\proclaim{Corollary 6}
Let $X$ be a superreflexive Banach space.  
With the same notations of the previous corollary, one has
$$
{\overline {\roman{span}}}^{\tau_u} 
\big\{ {\roman{Conv}}_b(X) \big\}
\supset
{\Cal{UC}}(X)
$$
(where $\tau_u$ is the topology of uniform convergence on $X$).
\endproclaim
\demo{Proof of the Corollary 6}
As we remarked during the proof of the previous {\bf Corollary 5}, 
the result is
proved if we show the $\tau_u$-density of the subset of 
$\Delta$-convex functions ${\roman{Conv}}_b(X)$
in the set of uniformly continuous functions $\Cal{UC}(X)$. 
\par
But this follows from {\bf Theorem 3} and Pisier's renorming theorem
(\cite {P}) that gives an equivalent norm with modulus of
convexity of power type $p$ (for some $p\geq 2$) 
on every superreflexive Banach space.
\remark{Remark}
For a simpler and more geometrical proof of Pisier's theorem 
we refer to \cite {L}.
\endremark
\enddemo
%
%
%%%%%%%%
\head 2.
The negative results
\endhead
%%%%%%%%

In this part, we will show that the  rotundity conditions on the norm 
needed in {\bf Theorem 1} can not be dropped. For instance,
even the pointwise convergence fails for some Banach spaces,
as the following counter-example shows.
\par
\example{Example 7}
There exists a Lipschitz function on $\ell_{\infty}$ which can not be a 
pointwise limit of a sequence of $\Delta$-convex functions.
\par
In the article \cite {T}, M. Talagrand proved that 
$
\text{Bor}( \ell_{\infty}, w )
\varsubsetneqq
\text{Bor}( \ell_{\infty} , {\|\cdot\|}_{\infty})
$.
Taking a $\|\cdot\|_{\infty}$-closed, non $w$-Borel set $B$ the function
$\text{d}(\cdot,B)$ is a Lipschitz function which can not be the pointwise
limit of any sequence of $\Delta$-convex functions (since $\Delta$-convex
functions are $w$-Borel). 
\endexample
\par
On the other hand, the $\tau_b$-density property of the span of 
${\text{Conv}}_b(X)$ in ${\Cal {UC}}_b$ 
is a characterization of superreflexivity. 
This conclusion comes from the following
theorem.
\par

\proclaim{Theorem 8}
For any non-superreflexive Banach space $X$ there exists a 
$1$--Lip\-schitz function defined on $X$ such that for every pair 
$\{c,d\}$ of
continuous, bounded on $B_X$, convex functions we have
$
\sup_{x\scriptstyle\in B_{X}} \big| f(x) - (c-d) (x) \big| \geq \frac {1}{4}
$.
\endproclaim
\demo{Proof of the Theorem 8}
Our main tool is James' Finite Tree Property (defined in \cite {J}). 
Specifically, we need the following well-known lemma.
\par
\proclaim{Lemma 8.1} 
Let $X$ be a non-superreflexive Banach space. 
For any  $0 < \theta < \frac {1}{2}$
and any $n\in\Bbb N$ there exists a dyadic $\roman{(}n,\theta\roman{)}$-tree
in the unit ball of $X$.
\endproclaim
\par
\demo{Proof of the lemma 8.1}
This lemma follows from the equivalence between the superreflexivity and
super--Radon-Nikodym properties (see \cite{Sn}). 
That means for any non-superreflexive Banach space $X$ 
that there exists a dual
Banach space failing the Radon-Nikodym Property which is finitely
representable on $X$ (namely, the bidual of one of its ultrapowers
${(X^{\Cal U})}^{**}$). The lemma is then deduced from
the existence of an infinite bounded dyadic tree in any dual Banach
space failing the Radon-Nikodym Property (\cite{Sg}).
\enddemo
This lemma allow us to construct a convenient sequence of trees in the 
unit ball of $X$. By ``convenient'' we understand the following:
\proclaim{Claim 8.2}
Let $\{\rho_n\}_{n\in\Bbb N}$ be a increasing sequence of positive numbers
tending to $1$. Then there exists a sequence of dyadic trees $T_{n}$ 
($n\in\Bbb N$) in the unit ball of $X$ such that
\roster
\item $T_{n}=\Big\{ x_{\alpha}^n : \alpha \in {\{-1,1\}}^{<n} \Big\}$ 
is a $\roman{(}n,\frac {{\rho}_{n}}{2} \roman{)}$-tree.
\item $\roman{dist}\roman{(}T_p, T_q \roman{)} \geq 
\frac {1}{2} {\rho}_{\max\{p,q\}}$.
\endroster
\endproclaim
Once the sequence ${\{T_n\}}_n$ of the previous claim is constructed, 
we define
the $1$-Lipschitz function $f(\cdot)=\text{dist}(\cdot, \cup_n S_n)$ where
$S_n=\{ x_\alpha \in T_n: | \alpha |  \text{\  is even} \}$.
\par
Obviously, if $x\in S_n$ then $f(x)=0$. On the other hand, the special
construction of the sequence $\{T_n\}$ gives that 
$f(x)\geq \frac {\rho_n}{2}$, for $x \in T_n \setminus  S_n$. 
Let us check that $f$
is the function we are looking for.
\par
If not, there is a pair of continuous convex functions $c,d$ bounded
on the unit ball of $X$ such that $|f(x)-(c-d)(x)|<\delta<\frac {1}{4}$,
$x\in B_X$. Take $n$ big enough so that $\delta<\frac {\rho_n}{4}$. Now,
we proceed to show that $f$ is strictly increasing along one of the branches
of $T_n=\{x_\alpha^n\}$ by a two-step algorithmic method.
\par
First, as $f(x_\varnothing^n)=0$ we have 
$d(x_\varnothing^n)\geq c(x_\varnothing^n)- \delta$. 
Then, 
using the convexity of $d$, we can find $\alpha_1 \in \{-1,1\}$ so that
$$
d(x_{\alpha_1}^n)\geq d(x_\varnothing^n) 
\geq c(x_\varnothing^n)- \delta. 
\tag{27}
$$ 
But, $\big( c(x_{\alpha_1}^n) - d(x_{\alpha_1}^n) \big)$ 
is $\delta$-close to
$f(x_{\alpha_1}^n)\geq \frac {\rho_n}{2}$, so we conclude from {\bf (27)} 
that 
$$
c(x_{\alpha_1}^n)
\geq 
f(x_{\alpha_1}^n) +  d(x_{\alpha_1}^n) - \delta
\geq
c(x_\varnothing^n) +  
\frac {\rho_n}{2} - 2\delta > 0.
\tag{28}
$$
Secondly, appealing this time to the convexity of $c$, we can choose 
$\alpha_2 \in \{-1,1\}$ such that 
$
c(x_{\alpha_1\roundhat\alpha_2}^n) 
\geq 
c(x_{\alpha_1}^n)
$.
Using {\bf (28)}, we deduce that
$$
c(x_{\alpha_1\roundhat\alpha_2}^n) - c(x_\varnothing^n) 
\geq 
\frac {\rho_n}{2} - 2\delta > 0.
\tag{28}
$$
And at this level we can repeat the same process as above.
\par
Iterating this process $n$ times up to the end of a branch of the tree
$T_n$, we find
a point 
$x_{\alpha_1\roundhat\cdots\roundhat\alpha_n}^n\in T_n$ 
that satisfies
$$
\multline
c(x_{\alpha_1\roundhat\cdots\roundhat\alpha_n}^n)-c(x_\varnothing^n)
\geq 
c(x_{\alpha_1\roundhat\cdots\roundhat\alpha_n}^n)
-
c(x_{\alpha_1\roundhat\cdots\roundhat\alpha_{n-2}}^n)
+ \cdots +
\\
+ \cdots +
c(x_{\alpha_1\roundhat\alpha_2}^n)-c(x_\varnothing^n)
\geq
\Big(\frac {\rho_n}{2} - 2\delta \Big) \frac {n}{2}.
\endmultline
$$ 
As $n$ can be chosen arbitrarily large, 
this last inequality contradicts the boundedness of $c$ in $B_X$.
\qed
\enddemo
\demo{Proof of the claim 8.2}
First, let us show that we can define a sequence 
${\{T_n\}}_{n\in \Bbb N}$ of trees 
such that for $n\in\Bbb N$ and the quotient map 
$\Pi_n:X\to X / {\text{span}\{T_1,\dots,T_n\}}$ 
one gets that $\Pi_n(T_{n+1})$ is a dyadic 
$(n+1,\frac {\rho_{n+1}}{2})$-tree rooted at $\overline{0}$.
By induction, consider that $T_1,\dots,T_n$ are already defined. 
Since $X$ is non-superreflexive, the quotient space 
$X / {\text{span}\{T_1,\dots,T_n\}}$
is non-superreflexive. So, there exists a dyadic $\overline{0}$-rooted
$(n+1,\frac {\rho_{n+1}}{2})$-tree 
$\overline{T_{n+1}}$ in  
$X / {\text{span}\{T_1,\dots,T_n\}}$. Then the tree $T_{n+1}$ which
we are looking for is simply  
a {\sl lifting} of the tree $\overline{T_{n+1}}$, obtained as follows: 
take 
$\big\{x_{\alpha}^{n+1}:\alpha \in {\{-1,1\}}^{n+1} \big\}\subset B_{X}$ 
in such a way that 
$\big\{ \Pi_{n}(x_{\alpha}^{n+1}) :  \alpha \in {\{-1,1\}}^{n+1} \big\}$
is the set of end points of $\overline{T_{n+1}}$;
then reconstruct from the ${\{x_{\alpha}^{n+1}\}}_{|\alpha|=n+1}$ 
the tree $T_{n+1}$ as 
$$
T_{n+1}:=
\Big\{ x_{\beta}^{n+1} = \sum\limits_{\beta \preceq \alpha}
 \frac {x_{\alpha}^{n+1}}{2^{n+1-|\beta|}} 
: \beta \in {\{-1,1\}}^{<n+1} \Big\}.
$$
\par
This sequence ${\{T_n\}}_{n\in \Bbb N}$ does not satisfy the conditions
of the claim yet; for example, the set of root points 
${ \{x_{\varnothing}^n \} }_{n}$ of the 
${\{T_n\}}_{n\in \Bbb N}$ might have cluster points. 
However, we can avoid this problem by selecting sub-trees of the trees
$T_n$ which do not contain these root points.
For each $n\in\Bbb N$, consider 
$$
T'_n=\{ y_\alpha^n : y_\alpha^{n}=x_{1\roundhat\alpha}^{n+1}\in T_{n+1} \}
$$
Let us see  
that the special quotient properties of $T_n$'s 
imply that if $p\neq q$ then 
$\text{dist}(T'_p, T'_q) \geq \frac {1}{2}{\rho}_{\max\{p,q\}}$.
For $q>p$ positive integers take $y_{\alpha}^{q} \in T'_q$ and 
$y_{\beta}^{p} \in T'_p$. Then we have that
$$
\Pi_{q}(y_{\beta}^{p})
=
\Pi_{q}(x_{1\roundhat\beta}^{p+1})
=
\overline{0}
=
\Pi_{q}(x_\varnothing^{q+1})
\text { and }
\Pi_{q}(y_{\alpha}^{q})
=
\Pi_{q}(x_{1\roundhat\alpha}^{q+1}).
\tag{30}
$$
Since $\Pi_{q}(T_{q+1})=\overline{T_{q+1}}$ is a dyadic 
$(q+1,\frac {\rho_{q+1}}2)$-tree rooted at $\overline{0}$ in 
the quotient space
$X/\text{span}\{T_1,\dots,T_{q-1}\}$, we deduce from {\bf (30)} that
$$
\| y_{\alpha}^{q} - y_{\beta}^{p} \| 
\geq
\|\Pi_q (y_{\alpha}^{q} - y_{\beta}^{p})\|
\geq
\| \Pi_{q}(x_\varnothing^{q+1}) - \Pi_{q}(x_{1\roundhat\alpha}^{q+1}) \|
\geq
\frac {\rho_{q+1}}2
\geq
\frac {1}{2}{\rho}_{\max\{p,q\}}.
\qed
$$ 
\enddemo
\par
%
%%%%%%%%%%%%%%%%%%%%%%%%%%%%%%%%%%%%%%%%%%%%%%%%%%%%%%%%%%%%%%%%%%%%%%%%%
%\newpage
%%%%%%%%%%%%%%%%%%%%%%%%%%%%%%%%%%%%%%%%%%%%%%%%%%%%%%%%%%%%%%%%%%%%%%%%%
%
\remark\nofrills{{\smc Acknowledgments.}}
\quad
The author would like to thank his advisor Professor Gilles Godefroy for 
many useful conversations.
The author is also indebted
to the Department of Mathematics 
of the University of Missouri-Columbia (U.S.A.)
for its hospitality during the last part of this work. 
\endremark
\refstyle{A}


\widestnumber\key{DFH$_1$}
\Refs

\ref\key{}
\by 
\paper
\jour
\vol
\yr
\issue
\pages
%\endref 
\book
\bookinfo
\publ
\yr
\publaddr
\endref

\ref\key{A}
\by E. Asplund
\paper Averaged norms
\jour Israel J. Math.
\vol 5
\yr 1967
\issue
\pages 227--233
\endref

\ref\key{DFH$_1$}
\by R. Deville-V. Fonf-P. H\'ajek
\paper Analytic and $C^k$ approximations of norms in separable Banach
spaces
\jour Studia Math.
\vol 120 
\yr 1996
\issue 1
\pages 61--74
\endref


\ref\key{DFH$_2$}
\bysame
\paper Analytic and polyhedral approximation of convex bodies in 
separable polyhedral spaces
\jour to appear
\vol
\yr
\issue
\pages
\endref 

\ref\key{DGZ}
\by R. Deville-G. Godefroy-V. Zizler
\book Smoothness and renormings in Banach spaces
\bookinfo Pitman Mono. and Surv. in Pure and App. Math.
\publ Longman
\yr 1993
\publaddr Boston
\endref

\ref\key{E}
\by G.A. Edgar
\paper Measurability in a Banach space
\jour Indiana Math. J.
\vol 26
\yr 1977
\issue 4
\pages 663-677
\endref

\ref\key{H}
\by J. Hoffman-J{\o}rgensen
\book On the Modulus of Smoothness and the G$_*$-Conditions in B-spaces
\bookinfo Preprint series
\publ Aarhus Univeritet, Matematisk Inst.
\yr 1974
\publaddr
\endref

\ref\key{J}
\by R.C. James
\book Some self-dual properties of normed linear spaces
\pages 159--175
\bookinfo  Symposium on Infinite-Dimensional Topology,
 Ann. of Math. Studies {\bf 69}
\publ Princenton Univ. Press.
\yr 1972
\publaddr Princenton
\endref

\ref\key{K}
\by J. Kurzweil 
\paper On approximation in real Banach spaces 
\jour Studia Math.
\vol 14
\yr 1954
\issue
\pages 214--231
\endref 

\ref\key{L}
\by G. Lancien
\paper On uniformly convex and uniformly Kadec-Klee renormings
\jour Serdica Math. J.
\vol 21
\yr 1995
\issue 1
\pages 1--18
\endref 

\ref\key{LL}
\by J.M. Lasry-P.L. Lions
\paper A remark on regularization in Hilbert spaces
\jour Israel J. Math.
\vol 55
\yr 1986
\pages 257--266
\endref

\ref\key{OS}
\by E. Odell-T. Schlumprecht
\paper The distortion problem
\jour Acta Math.
\vol 173
\yr 1994
\issue
\pages 259--281
\endref

\ref\key{Ph}
\by R.R. Phelps
\book Convex functions, Monotone operators and differentiability
\bookinfo Lecture Notes in Math. \vol 1364
\publ Springer-Verlag
\yr 1993
\publaddr Berlin
\endref

\ref\key{Pi}
\by G. Pisier
\paper Martingales with values in uniformly convex spaces
\jour Israel J. Math.
\vol 20
\yr 1975
\issue 
\pages 236--350
\endref 


\ref\key{PVZ}
\by R. Poliquin-J. Vanderwerff-V. Zizler
\paper Convex composite representation of lower semicontinuous functions
and renormings
\jour C. R. Acad. Sci. Paris, S\'erie I,
\vol 317
\yr 1993
\issue 
\pages 545--549
\endref

\ref\key{T}
\by     M. Talagrand
\paper  Comparaison des bor\'eliens d'un espace de Banach pour les topologies
        fortes et faibles 
\jour   Indiana Math. J.
\vol    27 
\yr   1978
\issue 6
\pages  1001--1004
\endref

\ref\key{Sg}
\by C. Stegall
\paper
The Radon-Nikodym property in conjugate Banach spaces
\jour Trans. Amer. Math. Soc.
\vol 206
\yr 1975
\pages 213--223
\endref

\ref\key{Sn}
\by J. Stern
\paper Propri\'et\'es locales et ultrapuissances d'espaces de Banach.
Expos\'es 7 et 8
\jour S\'emi\-naire Maurey-Schwartz (1974-1975). Centre Math.- \'Ecole
Polytech.
\year 1975
\publaddr Paris
\endref

\ref\key{St}
\by T. Str\"omberg
\paper On regularization in Banach spaces
\jour Ark. Mat.
\vol 34
\yr 1996
\issue
\pages 383--406
\endref

\endRefs
\enddocument